\documentclass[12pt,a4paper]{amsart}

\usepackage[utf8]{inputenc}
\usepackage[T1]{fontenc}

\usepackage[english]{babel}

\usepackage{amsmath}
\usepackage{amssymb}

\usepackage{cite}

\usepackage{enumerate}

\usepackage{12many}

\usepackage{pgfplots}
\pgfplotsset{compat=1.14}
\pgfplotscreateplotcyclelist{mylist}{{color=blue}, {color=green},
  {color=yellow}, {color=orange}, {color=red}}

\usepackage{tikz}
\usepackage{tikz-cd}

\title{Static Perfect Fluids with Symmetries}

\author{Adriano, L. $^{1}$} \address{ $^{1}$ Universidade Federal de
  Goiás, IME, 131, 74001-970, Goiânia, GO, Brazil.  }
\email{levi@ufg.br $^{1}$}

\author{Barboza, M. $^{2}$} \address{ $^{2}$ Instituto Federal Goiano,
  Rodovia Geraldo Silva Nascimento Km 2,5, 75790-000, Urutaí, GO,
  Brazil.  } \email{marcelo.barboza@ifgoiano.edu.br $^{2}$}

\author{Tokura, W. $^{3}$} \address{ $^{3}$ Universidade Federal de
  Goiás, IME, 131, 74001-970, Goiânia, GO, Brazil.  }
\email{williamisaotokura@hotmail.com $^{3}$}

\thanks{$^{3}$ Supported by CAPES/Brazil}

\date{\today}

\newcommand{\reals}{\mathbb{R}}
\newcommand{\positivereals}{(0,\infty)}

\newcommand{\hyperbolicspace}{\mathbb{H}^n}
\newcommand{\euclideanspace}{\mathbb{R}^n}
\newcommand{\sphericalspace}{\mathbb{S}^n}

\newcommand{\scalarfunctions}[1]{C^\infty(#1)}
\newcommand{\vectorfields}[1]{\mathfrak{X}^\infty(#1)}

\newcommand{\into}{\longrightarrow}
\newcommand{\goesto}{\longmapsto}

\newcommand{\dx}[1]{\frac{\partial}{\partial x_{#1}}}

\newcommand{\dxi}[1]{\xi_{,#1}}
\newcommand{\ddxi}[2]{\xi_{,#1#2}}

\newtheorem{example}{Example}
\newtheorem{proposition}{Proposition}
\newtheorem{remark}{Remark}
\newtheorem{theorem}{Theorem}
\newenvironment{myproof}[1]{\paragraph{\textbf{Proof of {#1}}}}{\hfill$\square$}

\DeclareMathOperator{\sech}{sech}

\begin{document}

\maketitle

\begin{abstract}
  In this paper we utilize symmetries in order to exhibit exact
  solutions to Einstein's equation of a perfect fluid on a static
  manifold whose spatial factor is conformal to a Riemannian space of
  constant sectional curvature. It's virtually possible to obtain
  infinitely many solutions via this approach since the equation is
  reduced into an ordinary differential equation that, essentially, is
  of the Riccati type. Three examples are shown in detail.
\end{abstract}


\section{Main Results}
\label{sec:intro}

Einstein's gravitational tensor of a Lorentzian manifold
$(\bar{M},\bar{g})$ is
\begin{equation}
  \label{eq:einsteins-tensor}
  G_{\bar{g}}=ric_{\bar{g}}-\frac{scal_{\bar{g}}}{2}\bar{g},
\end{equation}
where $ric_{\bar{g}}$ and $scal_{\bar{g}}$ stand for, respectively,
the Ricci tensor and the scalar curvature of the metric $\bar{g}$. It
should be noticed that the above tensor has zero divergence. General
relativity flows from Einstein's equation
\begin{equation}
  \label{eq:einsteins-equation}
  G_{\bar{g}}=T,
\end{equation}
of a manifold $(\bar{M},\bar{g})$ filled up with matter represented by
the stress-energy tensor $T$. According to general relativity, the
geometric properties of the universe are not independent, but rather
determind by matter. Therefore, it's only possible to infer something
about the geometric nature of the universe when the state of matter is
supposed to be known. We consider Einstein's equation of a perfect
fluid (see \cite{choquet-bruhat_general_2009},
\cite{hano_conformally-flatness_1981},
\cite{oneill_semi-riemannian_1983}) on $(\bar{M},\bar{g})$, i.e.,
\begin{equation}
  \label{eq:perfect-fluid}
  ric_{\bar{g}}-\frac{scal_{\bar{g}}}{2}\bar{g}=(\mu+\nu)\bar{g}(\cdot,X)\otimes\bar{g}(\cdot,X)+\nu\bar{g},
\end{equation}
where functions $\mu,\nu\in\scalarfunctions{\bar{M}}$ measure each a
specific feature of the fluid, those being energy density in the case
of $\mu$ and pressure in that of $\nu$, and the vector field
$X\in\vectorfields{\bar{M}}$, whose flow represents the dynamics of
the fluid, accomplishes \[\bar{g}(X,X)=-1,\] all along $\bar{M}$. The
definition of a perfect fluid, however, does not tell how to build out
a model of one. That is the reason why we have chosen to stick with
the so called static manifolds. A Lorentzian manifold
$(\bar{M},\bar{g})$ is called (globally) static (see
\cite{choquet-bruhat_general_2009},
\cite{hano_conformally-flatness_1981}) if
\begin{equation}
  \label{eq:static-manifold}
  \left\{
    \begin{aligned}
      \bar{M}&=M\times\reals, \\
      \bar{g}&=x^*g-f(x)^2dt^2,
    \end{aligned}
  \right.
\end{equation}
where
\[
  \begin{array}{ccccccccccc}
    x&:&\bar{M}&\into&M&\quad\mbox{and}\quad&t&:&\bar{M}&\into&\reals, \\
     &&(x,t)&\goesto&x&&&&(x,t)&\goesto&t
  \end{array}
\]
are the natural projections, $(M,g)$ is Riemannian and
$f\in\scalarfunctions{M}$ is positive.  It is known (see
\cite{hano_conformally-flatness_1981},
\cite{oneill_semi-riemannian_1983}) that under
\eqref{eq:static-manifold}, \eqref{eq:perfect-fluid} is equivalent to
\begin{equation}
  \label{eq:traceless}
  ric_g-\frac{scal_g}{n}g=\frac{1}{f}\Big((\nabla^2f)_g-\frac{(\Delta f)_g}{n}g\Big)\quad\mbox{on}\quad M,
\end{equation}
with relations
\begin{equation}
  \label{eq:energy-and-pressure}
  \mu=\frac{scal_g}{2}\quad\mbox{and}\quad\nu=\frac{n-1}{n}\frac{(\Delta f)_g}{f}-\frac{n-2}{n}\mu,
\end{equation}
where $(\nabla^2f)_g$ is the Hessian and $(\Delta f)_g$ is the
Laplacian of $f$ with respect to the metric $g$. In the scope of
conformal geometry we suppose that
\begin{equation}
  \label{eq:conformal-metric}
  (M,g)=(M^n_{\kappa},h^{-2}g_{\kappa}),
\end{equation}
where $(M^n_{\kappa},g_{\kappa})$ is a geodesically complete, simply
connected Riemannian space of constant sectional curvature
$\kappa\in\{-1,0,1\}$ and $h\in\scalarfunctions{M^n_{\kappa}}$ is a
positive function yet to be found. Therefore, we have that:
\begin{enumerate}[1.]
\item $(M^n_{0},g_{0})$ is the Euclidean $n$-space:
  \[M^n_{0}=\{(x_1,\ldots,x_n):
    x_1,\ldots,x_n\in\reals\}=\euclideanspace,\] with
  \[
    g_{0}\left(\dx{i},\dx{j}\right)=\delta_{ij}=\left\{
      \begin{array}{ccc}
        1&\mbox{if}&i=j,\\
        0&\mbox{if}&i\neq j,
      \end{array}
    \right.
  \]
  in the parametrization
  \[x:\euclideanspace\into\euclideanspace,\quad p\goesto p;\]
\item $(M^n_{-1},g_{-1})$ is the hyperbolic $n$-space:
  \[M^n_{-1}=\{(x_1,\ldots,x_n)\in\euclideanspace:
    x_{n}>0\}=\hyperbolicspace,\] with
  \[g_{-1}\left(\dx{i},\dx{j}\right)=x_{n}^{-2}\delta_{ij},\] in the
  parametrization
  \[x:\hyperbolicspace\into\hyperbolicspace,\quad p\goesto p,\] which
  is very similar in nature to that of the previous case;
\item $(M^n_{1},g_{1})$ is the euclidean $n$-sphere:
  \[M^n_{1}=\left\{(x_{1},\ldots,x_{n},x_{n+1})\in\reals^{n+1}\,:\,\sum\limits_{i=1}^{n+1}x_{i}^{2}=1\right\}=\sphericalspace,\]
  with
  \[g_{1}\left(\dx{i},\dx{j}\right)=\left(\frac{1+r}{2}\right)^{-2}\delta_{ij},\]
  in the parametrization given as the inverse of the stereographic
  projection with respect to either of the poles
  $\varepsilon e_{n+1}\in\sphericalspace$, $\varepsilon=\pm1$:
  \[x:\euclideanspace\into\sphericalspace\setminus\{-\varepsilon
    e_{n+1}\},\quad
    p\goesto\frac{2}{1+r}p+\varepsilon\frac{1-r}{1+r}e_{n+1},\] where
  $\{e_1,\ldots,e_n,e_{n+1}\}$ is the canonical linear basis of
  $\reals^{n+1}$,
  \[\euclideanspace\equiv\{(x_1,\ldots,x_n,x_{n+1})\in\reals^{n+1}:
    x_{n+1}=0\},\] and
  \[r:\euclideanspace\into[0,\infty),\quad
    (p_1,\ldots,p_n,0)\goesto\sum\limits_{i=1}^np_i^2.\]
\end{enumerate}
In order to handle all these geometries at once we let
$x:U\subset\euclideanspace\into M^n_{\kappa}$ represent the selected
parametrization according to the value of $\kappa\in\{-1,0,1\}$ and,
mostly important, we write
\[
  g_{\kappa}\left(\dx{i},\dx{j}\right)=\rho_\kappa^{-2}\delta_{ij},\]
where
\[
  \rho_\kappa(x_1,\ldots,x_n)=\left\{
    \begin{array}{lllll}
      1&\mbox{on}&\euclideanspace&\mbox{if}&\kappa=0,\\
      x_n&\mbox{on}&\hyperbolicspace&\mbox{if}&\kappa=-1,\\
      \frac{1+r}{2}&\mbox{on}&\euclideanspace&\mbox{if}&\kappa=1.\\
    \end{array}
  \right.
\]
Here as in \cite{barboza_invariant_2018} we assume that there do exist
functions
\begin{equation}
  \label{eq:invariant-functions}
  \xi:M^n_{\kappa}\into(a,b)\subset\reals\quad\mbox{and}\quad f,h:(a,b)\subset\reals\into\positivereals,
\end{equation}
making commutative diagrams out of those drawn below:
\begin{equation*}
  \begin{tikzcd}[row sep=huge, column sep=huge]
    M^n_{\kappa}\arrow[d,dashed,swap,"\xi"]\arrow[dr,"f\,\equiv\,f\,\circ\,\xi"] & & M^n_{\kappa}\arrow[d,dashed,swap,"\xi"]\arrow[dr,"h\,\equiv\,h\,\circ\,\xi"] & \\
    (a,b)\arrow[r,dashed,swap,"f"] & \positivereals &
    (a,b)\arrow[r,dashed,swap,"h"] & \positivereals
  \end{tikzcd}
\end{equation*}
We then get, by the chain rule, expressions
\[f_{,i}=\frac{\partial f}{\partial
    x_i}=\frac{df}{d\xi}\frac{\partial\xi}{\partial
    x_i}=f^{\prime}\xi_{,i}\quad\mbox{and}\quad
  f_{,ij}=\frac{\partial^2f}{\partial x_j\partial x_i}=\frac{\partial
    f_{,i}}{\partial
    x_j}=f^{\prime\prime}\xi_{,i}\xi_{,j}+f^{\prime}\xi_{,ij},\]
regarding the partial derivatives of the function $f$ up to the second
order (with analogous ones valid for $h$) in the parametrization
$x:U\subset\euclideanspace\into M^n_{\kappa}$. Now that everything has
been settled down, we believe it's time to state the main result of
this paper.

\begin{theorem}
  \label{thm:main}
  On each point of the parameter domain $U\subset\euclideanspace$, if
  it's not only true that
  \[\xi_{,ij}-\frac{(\Delta\xi)_{g_{0}}}{n}\delta_{ij}+\xi_{,i}\frac{\rho_{\kappa,j}}{\rho_\kappa}+\frac{\rho_{\kappa,i}}{\rho_\kappa}\xi_{,j}-\frac{2}{n}g_{0}\left((\nabla\xi)_{g_{0}},\frac{(\nabla{\rho_\kappa})_{g_{0}}}{\rho_\kappa}\right)\delta_{ij}=0,\]
  for all $i,j\in\{1,\ldots,n\}$, but also that
  \[\xi_{,k}\xi_{,l}-\frac{\|(\nabla\xi)_{g_{0}}\|^2}{n}\delta_{kl}\neq0,\]
  for some $k,l\in\{1,\dots,n\}$, then $f$ is a solution to
  \eqref{eq:traceless} on $(M,g)=(M^n_{\kappa},h^{-2}g_{\kappa})$ if,
  and only if,
  \[(n-2)\frac{h^{\prime\prime}}{h}-2\frac{h^{\prime}}{h}\frac{f^{\prime}}{f}-\frac{f^{\prime\prime}}{f}=0,\]
  on $\xi(x(U))=\{\xi(q)\mid q\in x(U)\}\subset(a,b)$. Energy density
  and pressure of the fluid are then given by:
  \[\mu=\frac{n-1}{2}\left(\kappa nh^2+2h(\Delta
      h)_{g_\kappa}-n\|(\nabla h)_{g_\kappa}\|^2\right),\] and
  \[\nu=\frac{n-1}{n}\frac{(\Delta
      f)_{g_\kappa}}{f}-\frac{n-2}{n}\mu,\] respectively.
\end{theorem}

\begin{remark}
  Upon declaring
  \[x=\frac{h^\prime}{h}\quad\mbox{and}\quad y=\frac{f^\prime}{f},\]
  we might observe that
  \begin{equation}
    \label{eq:ode}
    (n-2)\frac{h^{\prime\prime}}{h}-2\frac{h^\prime}{h}\frac{f^\prime}{f}-\frac{f^{\prime\prime}}{f}=0,
  \end{equation}
  becomes
  \begin{equation}
    \label{eq:riccati}
    y^\prime=(n-2)(x^2+x^\prime)-2xy-y^2,
  \end{equation}
  which is Riccati in $y$ for a known $x$. Therefore, it's general
  solution is of the form $y=y_0+u$, where $y_0$ is a particular
  solution of \eqref{eq:riccati} and $u$ must solve the linear
  equation
  \[\frac{d}{d\xi}\left(\frac{1}{u}\right)-2(y_0+x)\frac{1}{u}=1.\]
\end{remark}

\begin{proposition}
  \label{prop:invariants}
  If the local expression of $\xi:M^n_\kappa\into\reals$ with respect
  to the parametrization $x:U\subset\euclideanspace\into M^n_\kappa$
  satisfies
  \[\xi_{,ij}-\frac{(\Delta\xi)_{g_{0}}}{n}\delta_{ij}+\xi_{,i}\frac{\rho_{\kappa,j}}{\rho_\kappa}+\frac{\rho_{\kappa,i}}{\rho_\kappa}\xi_{,j}-\frac{2}{n}g_{0}\left((\nabla\xi)_{g_{0}},\frac{(\nabla{\rho_\kappa})_{g_{0}}}{\rho_\kappa}\right)\delta_{ij}=0,\]
  on each point of the parameter domain $U\subset\euclideanspace$ for
  all $i,j\in\nto{1}{n}$, then:
  \begin{enumerate}[1.]
  \item $\kappa=0$:
    \[\xi:\euclideanspace\into\reals,\quad(x_1,\ldots,x_n)\goesto\sum\limits_{i=1}^{n}\left(\frac{a}{2}x_i^2+b_ix_i+c_i\right),\]
    where $a,b_1,\ldots,b_n,c_1,\ldots,c_n\in\reals$ are constants;
  \item $\kappa=-1$:
    \[\xi:\hyperbolicspace\into\reals,(x_1,\ldots,x_n)\goesto\frac{1}{x_n}\sum\limits_{i=1}^{n}\left(\frac{a}{2}x_i^2+b_ix_i+c_i\right),\]
    where $a,b_1,\ldots,b_n,c_1,\ldots,c_n\in\reals$ are constants;
  \item $\kappa=1$:
    \[\xi:\sphericalspace\into\reals,\quad
      (p_1,\ldots,p_n,p_{n+1})\goesto p_{n+1},\]
  \end{enumerate}
  if, in addition, it's assumed that there do exist
  \[a<0\quad\mbox{and}\quad\xi:(a,\infty)\subset\reals\into\reals,\]
  such that
  \[\xi(x)=\xi(r(x)),\]
  for all $x\in\euclideanspace$, where
  \[r:\euclideanspace\into[0,\infty),\quad(x_1,\ldots,x_n)\goesto\sum\limits_{i=1}^{n}x_i^2.\]
\end{proposition}


\section{Examples}
\label{sec:examples}

\begin{example}
  As for the case $\kappa=0$, we might choose
  \[\xi:\euclideanspace\into\reals,\quad (p_1,\ldots,p_n)\goesto
    \frac{1}{\sqrt{n}}p_1+\frac{1}{\sqrt{n}}p_2+\cdots+\frac{1}{\sqrt{n}}p_n,\]
  and then
  \[h:\euclideanspace\into\reals,\quad p\goesto\cos(\xi(p)),\] from
  which we see that $\euclideanspace$ must be shrinked into one of its
  open subsets, let's say,
  \[M^n=\left\{p\in\euclideanspace\,:\,\,-\frac{\pi}{2}<\xi(p)<\frac{\pi}{2}\right\},\]
  so that $h$ gets to be strictly positive. This leaves us with the
  equation
  \[f^{\prime\prime}-2\tan(\xi)f^\prime+(n-2)f=0,\] of which
  \[f:M^n\into\positivereals,\quad
    p\goesto\frac{e^{-\xi(p)\sqrt{n-1}}}{\cos(\xi(p))},\] is a
  positive solution. Therefore, the manifold
  \[\bar{M}=M\times\reals,\]
  munished with the metric tensor
  \[\bar{g}=\sec^2(\xi(x))\left[x^*g_0-\left(e^{-\xi(x)\sqrt{n-1}}\right)^2dt^2\right],\]
  solves Einstein's equation for the perfect fluid
  \[
    \begin{aligned}
      T&=(\mu+\nu)\bar{g}(\cdot,X)\otimes\bar{g}(\cdot,X)+\nu\bar{g}\\
      &=\sec^2(\xi(x))\left[\nu(\xi(x))
        x^*g_0+\mu(\xi(x))\left(e^{-\xi(x)\sqrt{n-1}}\right)^2dt^2\right],
    \end{aligned}
  \]
  characterized by its energy density
  \[\mu(\xi)=\frac{n-1}{2}\left((n-2)\cos^2(\xi)-n\right),\]
  pressure
  \[
    \begin{aligned}
      \nu(\xi)&=\frac{(n-1)\cos^2(\xi)}{n}\Big\{[\tan(\xi)+\sqrt{n-1}]^2+\sec^2(\xi)+\\
      &+(n-2)\tan(\xi)[\tan(\xi)-\sqrt{n-1}]\Big\}-\frac{n-2}{n}\mu(\xi),
    \end{aligned}
  \]
  and the vector field $X=\frac{1}{f}\frac{\partial}{\partial t}$
  where, in all of the above, $x:\bar{M}\into M$ and
  $t:\bar{M}\into\reals$ denote the natural projections.
\end{example}

\begin{example}
  By choosing
  \[\xi:\hyperbolicspace\into\reals,\quad(x_1,\ldots,x_n)\goesto\frac{1}{x_n}\sum\limits_{i=1}^{n-1}b_ix_i,\]
  where $b_1,\ldots,b_{n-1}\in\reals$ are constants subject to
  \[\sum\limits_{i=1}^{n-1}b_i^2=1,\]
  and also
  \[h:\hyperbolicspace\into\positivereals,\quad
    x\goesto\cosh(\xi(x)),\] as representatives of the case
  $\kappa=-1$, we get the equation
  \[f^{\prime\prime}+2\tanh(\xi)f^\prime-(n-2)f=0,\] of which
  \[f:\hyperbolicspace\into\positivereals,\quad
    x\goesto\frac{e^{-\xi(x)\sqrt{n-1}}}{\cosh(\xi(x))},\] is a
  positive solution. Thus, the manifold
  \[\bar{M}=\hyperbolicspace\times\reals,\]
  equipped with the metric tensor
  \[\bar{g}=\sech(\xi(x))\left[x^*g_{-1}-(e^{-\xi(x)\sqrt{n-1}})^2dt^2\right],\]
  solves Einstein's equation for the perfect fluid
  \[
    \begin{aligned}
      T&=(\mu+\nu)\bar{g}(\cdot,X)\otimes\bar{g}(\cdot,X)+\nu\bar{g}\\
      &=\sech^2(\xi(x))\left[\nu(\xi(x))
        x^*g_{-1}+\mu(\xi(x))\left(e^{-\xi(x)\sqrt{n-1}}\right)^2dt^2\right],
    \end{aligned}
  \]
  characterized by its energy density
  \[
    \begin{aligned}
      \mu(\xi)&=\frac{n-1}{2}\Big\{-n\cosh^2(\xi)+2\cosh(\xi)[(1+\xi^2)\cosh(\xi)+n\xi\sinh(\xi)]+\\
      &-n(1+\xi^2)\sinh^2(\xi)\Big\},
    \end{aligned}
  \]
  pressure
  \[
    \begin{aligned}
      \nu(\xi)&=\frac{(n-1)\cosh^2(\xi)}{n}\Big\{(1+\xi^2)[(\tanh(\xi)+\sqrt{n-1})^2-\sech^2(\xi)]+\\
      &+n\xi(\tanh(\xi)+\sqrt{n-1})+(n-2)(1+\xi^2)\tanh(\xi)(\tanh(\xi)+\sqrt{n-1})\Big\}+\\
      &-\frac{n-2}{n}\mu(\xi)
    \end{aligned}
  \]
  and the vector field $X=\frac{1}{f}\frac{\partial}{\partial t}$
  where, as usual, $x:\bar{M}\into\hyperbolicspace$ and
  $t:\bar{M}\into\reals$ indicate the natural projections.
\end{example}

\begin{example}
  Given that we have
  \[\xi:\sphericalspace\into\reals,\quad
    (p_1,\ldots,p_n,p_{n+1})\goesto p_{n+1},\] in case $\kappa=1$, we
  see that upon choosing
  \[h:\sphericalspace\into\positivereals,\quad p\goesto\cos(\xi(p)),\]
  we are lead to the equation
  \[f^{\prime\prime}-2\tan(\xi)f^\prime+(n-2)f=0,\] of which
  \[f:\sphericalspace\into\positivereals,\quad p\goesto
    \frac{e^{-\xi(p)\sqrt{n-1}}}{\cos(\xi(p))},\] is a positive
  solution. Henceforth, the manifold
  \[\bar{M}=\sphericalspace\times\reals,\]
  furnished with metric tensor
  \[\bar{g}=\sec^2(\xi(x))\left[x^*g_1-\left(e^{-\xi(x)\sqrt{n-1}}\right)^2dt^2\right],\]
  solves Einstein's equation for the perfect fluid
  \[
    \begin{aligned}
      T&=(\mu+\nu)\bar{g}(\cdot,X)\otimes\bar{g}(\cdot,X)+\nu\bar{g}\\
      &=\sec^2(\xi(x))\left[\nu(\xi(x))
        x^*g_1+\mu(\xi(x))\left(e^{-\xi(x)\sqrt{n-1}}\right)^2dt^2\right],
    \end{aligned}
  \]
  characterized by its energy density
  \[
    \begin{aligned}
      \mu(\xi)&=\frac{(n-1)}{n}\Big\{n\cos^2(\xi)+2\cos(\xi)[-(1-\xi^2)\cos(\xi)+n\xi\sin(\xi)]+\\
      &-n(1-\xi^2)\sin^2(\xi)\Big\}
    \end{aligned}
  \]
  pressure
  \[
    \begin{aligned}
      \nu(\xi)&=\frac{(n-1)\cos^2(\xi)}{n}\Big\{(1-\xi^2)[(\tan(\xi)-\sqrt{n-1})^2+\sec^2(\xi)]+\\
      &-n\xi(\tan(\xi)-\sqrt{n-1})+(n-2)(1-\xi^2)\tan(\xi)(\tan(\xi)-\sqrt{n-1})\Big\}+\\
      &-\frac{n-2}{n}\mu(\xi),
    \end{aligned}
  \]
  and vector field $X=\frac{1}{f}\frac{\partial}{\partial t}$ where,
  once more, $x:\bar{M}\into\sphericalspace$ and
  $t:\bar{M}\into\reals$ stand for the natural projections. Next, we
  plot the graphs of both $\mu$ and $\nu$ as functions of $\xi$ in the
  dimensions $n=2,3,4,5$ and $6$.
  \[
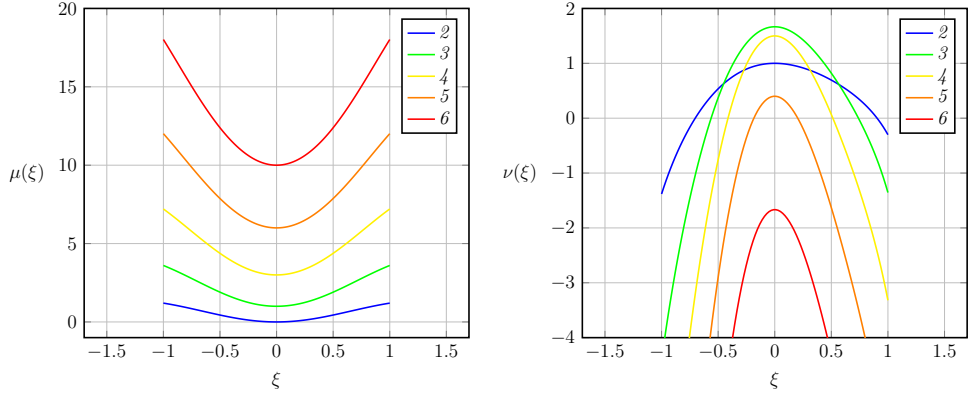
    \begin{array}{cc}
      \begin{tikzpicture}[scale=0.6]
        \begin{axis}[ grid=major , xmin=-1.7, xmax=1.7 , ymin=-1,
          ymax=20 , xlabel=$\xi$,
          ylabel=$\mu(\xi)$ , ylabel style={rotate=-90} , domain=-1:1
          , samples=100 , cycle list name=mylist , line width=1pt, ,
          scale only axis , legend style={legend pos=north east} ]
          \foreach \n in {2,...,6}{
            \addplot{(\n-1)*(\n*cos(deg(2*x))+\n*x*sin(deg(2*x))+\n*x^2-2*cos(deg(x))^2-(\n-2)*x^2*cos(deg(x))^2)/2};
            \expandafter\addlegendentry\expandafter{\n}; }
        \end{axis}
      \end{tikzpicture}
      &
        \begin{tikzpicture}[scale=0.6]
          \begin{axis}[ grid=major , xmin=-1.7, xmax=1.7, ymin=-4,
            ymax=2 , xlabel=$\xi$, ylabel=$\nu(\xi)$ , ylabel
            style={rotate=-90} , domain=-1:1 , samples=100 , cycle
            list name=mylist , line width=1pt, , scale only axis ,
            legend style={legend pos=north east} ]
            \foreach \n in {2,...,6}{
              \addplot{((\n-1)*cos(deg(x))*((1-x^2)*(\n-2*sqrt(\n-1)*tan(deg(x))+2*tan(deg(x))^2)-\n*x*(-sqrt(\n-1)+tan(deg(x)))+(\n-2)*(1-x^2)*tan(deg(x))*(-sqrt(\n-1)+tan(deg(x)))))/\n-(\n-2)*((\n-1)*(\n*cos(deg(2*x))+\n*x*sin(deg(2*x))+\n*x^2-2*cos(deg(x))^2-(\n-2)*x^2*cos(deg(x))^2)/2)/\n};
              \expandafter\addlegendentry\expandafter{\n}; }
          \end{axis}
        \end{tikzpicture}
    \end{array}
  \]
\end{example}


\section{Proofs}
\label{sec:proofs}

Please, take a look at section \ref{sec:intro} to familiarize yourself
with both the notation we adopt and the various conventions we do
make.

\vspace{12pt}

\begin{myproof}{Theorem \ref{thm:main}}
  Exactly how the change of metrics $g_\kappa\goesto g=h^{-2}g_\kappa$
  seems to alter the tensors $ric_{g_\kappa}=\kappa(n-1)g_\kappa$ and
  $(\nabla^2f)_{g_\kappa}$ is something better comprehended with the
  help of the next $2$ formulas (see \cite{besse_einstein_2008}):
  \begin{equation}
    \label{eq:conformal-ric}
    \begin{aligned}
      ric_g&=\kappa(n-1)g_\kappa+\\
      &+h^{-2}\Big\{(n-2)h(\nabla^2h)_{g_\kappa}+\Big[h(\Delta
      h)_{g_\kappa}-(n-1)\|(\nabla
      h)_{g_\kappa}\|^2\Big]g_\kappa\Big\},
    \end{aligned}
  \end{equation}
  and
  \begin{equation}
    \label{eq:conformal-hess}
    \begin{aligned}
      (\nabla^2f)_g&=(\nabla^2f)_{g_\kappa}+\\
      &+df\otimes d(\log h)+d(\log h)\otimes df-g_{\kappa}((\nabla
      f)_{g_\kappa},(\nabla\log h)_{g_\kappa}).
    \end{aligned}
  \end{equation}
  Therefore, it's readily seen that
  \begin{equation}
    \label{eq:conformal-scal}
    scal_g=\kappa n(n-1)h^2+(n-1)\Big[2h(\Delta h)_{g_\kappa}-n\|(\nabla h)_{g_\kappa}\|^2\Big],
  \end{equation}
  and
  \begin{equation}
    \label{eq:conformal-lapl}
    (\Delta f)_g=h^2\Big[(\Delta f)_{g_\kappa}-(n-2)g_\kappa((\nabla f)_{g_\kappa},(\nabla\log h)_{g_\kappa})\Big].
  \end{equation}
  From expressions \eqref{eq:conformal-ric} and
  \eqref{eq:conformal-scal} comes
  \begin{equation}
    \label{eq:conformal-traceless-ric}
    ric_g-\frac{scal_g}{n}g=\frac{n-2}{h}\left[(\nabla^2h)_{g_\kappa}-\frac{(\Delta h)_{g_\kappa}}{n}g_\kappa\right],
  \end{equation}
  whilst
  \begin{equation}
    \label{eq:conformal-traceless-hess}
    \begin{aligned}
      &(\nabla^2f)_g-\frac{(\Delta f)_g}{n}g=(\nabla^2 f)_{g_\kappa}-\frac{(\Delta f)_{g_\kappa}}{n}g_\kappa+\\
      &+df\otimes d(\log h)+d(\log h)\otimes
      df-\frac{2}{n}g_\kappa((\nabla f)_{g_\kappa},(\nabla\log
      h)_{g_\kappa})g_\kappa.
    \end{aligned}
  \end{equation}
  follows from \eqref{eq:conformal-hess} and
  \eqref{eq:conformal-lapl}. By \eqref{eq:conformal-traceless-hess},
  since
  \[g_\kappa=\rho_\kappa^{-2}g_0,\] on the open set
  $x(U)\subset M^n_\kappa$, which is the image set of the
  parametrization $x:U\subset\euclideanspace\into M^n_\kappa$,
  \eqref{eq:conformal-traceless-ric} becomes
  \begin{equation}
    \label{eq:conformal-traceless-ric-2}
    \begin{aligned}
      &ric_g-\frac{scal_g}{n}g=\frac{n-2}{h}\Bigg[(\nabla^2h)_{g_0}-\frac{(\Delta
        h)_{g_0}}{n}g_0+\\
      &+dh\otimes d(\log\rho_\kappa)+d(\log\rho_\kappa)\otimes
      dh-\frac{2}{n}g_0\left((\nabla
        h)_{g_0},(\nabla\log\rho_\kappa)_{g_0}\right)g_0\Bigg],
    \end{aligned}
  \end{equation}
  and \eqref{eq:conformal-traceless-hess} itself gives
  \begin{equation}
    \label{eq:conformal-traceless-hess-2}
    \begin{aligned}
      &(\nabla^2f)_g-\frac{(\Delta f)_g}{n}g=(\nabla^2
      f)_{g_0}-\frac{(\Delta f)_{g_0}}{n}g_0+\\
      &+df\otimes d(\log\rho_\kappa)+d(\log\rho_\kappa)\otimes
      df-\frac{2}{n}g_0((\nabla
      f)_{g_0},(\nabla\log\rho_\kappa)_{g_0})g_0+\\
      &+df\otimes d(\log h)+d(\log h)\otimes df-\frac{2}{n}g_0((\nabla
      f)_{g_0},(\nabla\log h)_{g_0})g_0,
    \end{aligned}
  \end{equation}
  where, in the last expression, we have used the fact that
  \[(\nabla f)_{g_\kappa}=\rho_\kappa^2(\nabla
    f)_{g_0}\quad\mbox{and}\quad(\nabla\log
    h)_{g_\kappa}=\rho_\kappa^2(\nabla\log h)_{g_0}.\] As a result,
  \eqref{eq:traceless} asserts that
  \begin{equation*}
    \left\{
      \begin{aligned}
        &\frac{n-2}{h}\Bigg[(\nabla^2h)_{g_0}-\frac{(\Delta h)_{g_0}}{n}g_0+\\
        &+dh\otimes d(\log\rho_\kappa)+d(\log\rho_\kappa)\otimes
        dh-\frac{2}{n}g_0\left((\nabla
          h)_{g_0},(\nabla\log\rho_\kappa)_{g_0}\right)g_0\Bigg]=\\
        &=\frac{1}{f}\Bigg[(\nabla^2f)_{g_0}-\frac{(\Delta f)_{g_0}}{n}g_0+\\
        &+df\otimes d(\log\rho_\kappa)+d(\log\rho_\kappa)\otimes
        df-\frac{2}{n}g_0((\nabla
        f)_{g_0},(\nabla\log\rho_\kappa)_{g_0})g_0+\\
        &+df\otimes d(\log h)+d(\log h)\otimes
        df-\frac{2}{n}g_0((\nabla f)_{g_0},(\nabla\log
        h)_{g_0})g_0\Bigg],
      \end{aligned}
    \right.
  \end{equation*}
  and because we have
  \begin{equation*}
    (\nabla f)_{g_0}=f^\prime(\nabla\xi)_{g_0}\quad\mbox{and}\quad(\Delta f)_{g_0}=f^{\prime\prime}\|(\nabla\xi)_{g_0}\|^2+f^\prime(\Delta\xi)_{g_0},
  \end{equation*}
  even that
  \begin{equation*}
    \begin{aligned}
      &\left(\xi_{,i}\xi_{,j}-\frac{\|(\nabla\xi)_{g_0}\|^2}{n}\delta_{ij}\right)\left[(n-2)\frac{h^{\prime\prime}}{h}-2\frac{h^\prime}{h}\frac{f^\prime}{f}-\frac{f^{\prime\prime}}{f}\right]=\\
      &=\left[\frac{f^\prime}{f}-(n-2)\frac{h^\prime}{h}\right]\left[\xi_{,ij}+\xi_{,i}\frac{\rho_{\kappa,j}}{\rho_\kappa}+\frac{\rho_{\kappa,i}}{\rho_\kappa}\xi_{,j}-\frac{2}{n}g_0\left((\nabla\xi)_{g_0},\frac{(\nabla\rho_\kappa)_{g_0}}{\rho_\kappa}\right)\right],
    \end{aligned}
  \end{equation*}
  for all $i,j\in\{1,\ldots,n\}$.
\end{myproof}

\vspace{10pt}

\begin{myproof}{Proposition \ref{prop:invariants}}\textup{}

  \vspace{5pt}\noindent\textbf{Case} $ \kappa=0 $: From the fact that
  \[
    \begin{aligned}
      0&=\ddxi{i}{j}-\frac{(\Delta\xi)_{g_0}}{n}\delta_{ij}+\dxi{i}\frac{\rho_{0,j}}{\rho_0}+\frac{\rho_{0,i}}{\rho_0}\dxi{j}-\frac{2}{n}g_0\left((\nabla\xi)_{g_0},\frac{(\nabla\rho_0)_{g_0}}{\rho_0}\right)\delta_{ij}\\
      &=\ddxi{i}{j}-\frac{(\Delta\xi)_{g_0}}{n}\delta_{ij},
    \end{aligned}
  \]
  for every $i,j\in\nto{1}{n}$, we get that $\ddxi{i}{j}=0$ whenever
  $i\neq j$. Thus,
  \[\xi(x_1,\ldots,x_n)=\sum\limits_{i=1}^nF_i(x_i),\]
  and, as such,
  \[F_i^{\prime\prime}(x_i)=\frac{1}{n}\sum\limits_{j=1}^nF_j^{\prime\prime}(x_j),\]
  for all $i\in\nto{1}{n}$. If we let
  \[F_1^{\prime\prime}=\cdots=F_n^{\prime\prime}=a\in\reals,\]
  then
  \[\xi(x_1,\ldots,x_n)=\sum\limits_{i=1}^n\left(\frac{a}{2}x_i^2+b_ix_i+c_i\right),\]
  where $b_1,\ldots,b_n,c_1,\ldots,c_n\in\reals$ are constants.
  
  \vspace{6pt}\noindent\textbf{Case} $ \kappa=-1 $: Notice that
  \[0=\ddxi{i}{j}-\frac{(\Delta\xi)_{g_{0}}}{n}\delta_{ij}+\dxi{i}\frac{\rho_{-1,j}}{\rho_{-1}}+\frac{\rho_{-1,i}}{\rho_{-1}}\dxi{j}-\frac{2}{n}g_0\left((\nabla\xi)_{g_0},\frac{(\nabla\rho_{-1})_{g_{0}}}{\rho_{-1}}\right)\delta_{ij},\]
  turns out to be
  \begin{equation}
    \label{eq:hyperbolic-invariants}
    0=\ddxi{i}{j}-\frac{(\Delta\xi)_{g_0}}{n}\delta_{ij}+\dxi{i}\frac{\delta_{jn}}{x_n}+\frac{\delta_{in}}{x_n}\dxi{j}-\frac{2}{n}\frac{\dxi{n}}{x_n}\delta_{ij},
  \end{equation}
  for all $i,j\in\nto{1}{n}$, now that $\kappa=-1$. Therefore,
  $1\leqslant i\leqslant n-1$ and $j=n$ together give
  \[(\xi_{,i}\cdot
    x_n)_{,n}=\left(\xi_{,in}+\frac{\xi_{,i}}{x_n}\right)x_n=0,\] from
  which we get that
  \[\dxi{i}(x_1,\ldots,x_n)=\frac{F_i(x_1,\ldots,x_{n-1})}{x_n},\]
  whenever $i\in\nto{1}{n-1}$. But then,
  \[\frac{F_{i,j}}{x_n}=\ddxi{i}{j}=0,\]
  for all $1\leqslant i\neq j\leqslant n-1$, culminating in
  $F_i=F_i(x_i)$ for every $i\in\nto{1}{n-1}$. Upon taking
  $1\leqslant i=j\leqslant n-1$ we see that
  \[\ddxi{i}{i}-\frac{(\Delta\xi)_{g_0}}{n}-\frac{2}{n}\frac{\dxi{n}}{x_n}=0,\]
  implying that
  \[\frac{F_i^\prime}{x_n}=\ddxi{i}{i}=\frac{(\Delta\xi)_{g_0}}{n}+\frac{2}{n}\frac{\dxi{n}}{x_n}=\ddxi{j}{j}=\frac{F_j^\prime}{x_n},\]
  for every $i,j\in\nto{1}{n-1}$. Henceforth, there must be an
  $a\in\reals$ such that
  \[F_1^\prime=\cdots=F_{n-1}^\prime=a,\] meaning that there does
  exist some $b_i\in\reals$ validating the identity
  \[F_i(x_i)=ax_i+b_i,\] for each $i\in\nto{1}{n-1}$. So far, we may
  solely guarantee that
  \[\xi(x_1,\ldots,x_n)=\frac{1}{x_n}\sum\limits_{i=1}^{n-1}\left(\frac{a}{2}x_i^2+b_ix_i+c_i\right)+G(x_n),\]
  where $c_1,\ldots,c_{n-1}\in\reals$ are constants but, as $i=j=n$
  lead to
  \[\ddxi{n}{n}-\frac{(\Delta\xi)_{g_0}}{n}+\frac{\dxi{n}}{x_n}+\frac{\dxi{n}}{x_n}-\frac{2}{n}\frac{\dxi{n}}{x_n}=0,\]
  once more by \eqref{eq:hyperbolic-invariants}, we have that
  \[\ddxi{n}{n}+2\frac{\dxi{n}}{x_n}=\frac{(\Delta\xi)_{g_0}}{n}+\frac{2}{n}\frac{\dxi{n}}{x_n}=\ddxi{i}{i}=\frac{a}{x_n},\]
  so that
  \[(\dxi{n}\cdot
    x_n^2)_{,n}=\left(\ddxi{n}{n}+2\frac{\dxi{n}}{x_n}\right)x_n^2=ax_n.\]
  It's now readily seen that
  \[\left[G^\prime\cdot
      x_n^2\right]_{,n}=\left[\left(-\frac{1}{x_n^2}\sum\limits_{i=1}^{n-1}\left(\frac{a}{2}x_i^2+b_ix_i\right)+G^\prime\right)\cdot
      x_n^2\right]_{,n}=ax_n,\] of
  which \[G(x_n)=\frac{a}{2}x_n+\frac{c_n}{x_n}+b_n,\] is the general
  solution, where $b_n,c_n\in\reals$ are constants. In conclusion, we
  have the function
  \[\xi(x_1,\ldots,x_n)=\frac{1}{x_n}\sum\limits_{i=1}^{n}\left(\frac{a}{2}x_i^2+b_ix_i+c_i\right).\]
  \textbf{Case} $ \kappa=1 $: As for the last case, there remains to
  solve the equation
  \[
    \begin{aligned}
      0&=\ddxi{i}{j}-\frac{(\Delta\xi)_{g_0}}{n}\delta_{ij}+\dxi{i}\frac{\rho_{1,j}}{\rho_1}+\frac{\rho_{1,i}}{\rho_1}\dxi{j}-\frac{2}{n}g_0\left((\nabla\xi)_{g_0},\frac{(\nabla\rho_1)_{g_0}}{\rho_1}\right)\delta_{ij}\\
      &=\ddxi{i}{j}-\frac{(\Delta\xi)_{g_0}}{n}\delta_{ij}+\dxi{i}\frac{2x_j}{1+r}+\frac{2x_i}{1+r}\dxi{j}-\frac{2}{n}g_0\left((\nabla\xi)_{g_0},\frac{2x}{1+r}\right)\delta_{ij},
    \end{aligned}
  \]
  and since it is a difficult one to get over with in all of its
  generality, we may tacitly assume that there does exist some
  function $\xi:(a,\infty)\subset\reals\into\positivereals$, where
  $a<0$, such that the local expression of
  $\xi:\sphericalspace\into\reals$ with respect to the parametrization
  $x:\euclideanspace\into\sphericalspace\setminus\{-\varepsilon
  e_{n+1}\}$ satisfies
  \[\xi(u)=\xi(r(u)),\] for all $u\in\euclideanspace$, where
  \[r:\euclideanspace\into\reals,\quad
    (u_1,\ldots,u_n)\goesto\sum\limits_{i=1}^{n}u_i^2.\] Because
  \[\dxi{i}=\xi^{\prime}2x_i\quad\mbox{and}\quad\ddxi{i}{j}=\xi^{\prime\prime}4x_ix_j+\xi^{\prime}2\delta_{ij},\]
  for every $i,j\in\nto{1}{n}$, we know that
  \[(\nabla\xi)_{g_0}=2\xi^\prime
    x\quad\mbox{and}\quad(\Delta\xi)_{g_0}=4r\xi^{\prime\prime}+2n\xi^\prime.\]
  Henceforth,
  \[
    \begin{aligned}
      0&=\xi^{\prime\prime}4x_ix_j+\xi^\prime2\delta_{ij}-\frac{\xi^{\prime\prime}4r+\xi^\prime2n}{n}\delta_{ij}+\xi^\prime\frac{4x_ix_j}{1+r}+\xi^\prime\frac{4x_ix_j}{1+r}-\frac{2}{n}\xi^\prime\frac{4r}{1+r}\delta_{ij}\\
      &=4\left(\xi^{\prime\prime}+\frac{2}{1+r}\xi^\prime\right)\left(x_ix_j-\frac{r}{n}\delta_{ij}\right).
    \end{aligned}
  \]
  Notice that
  \[\xi:\sphericalspace\into\reals,\quad
    (p_1,\ldots,p_n,p_{n+1})\goesto p_{n+1},\] owns the expression
  \[\xi(r)=\frac{2\varepsilon}{1+r}-\varepsilon,\]
  with respect to
  \[x:\euclideanspace\into\sphericalspace\setminus\{-\varepsilon
    e_{n+1}\},\quad u\goesto
    \frac{2}{1+r}u+\varepsilon\frac{1-r}{1+r}e_{n+1},\] with which
  it's seen that
  \[\xi^\prime=\frac{d\xi}{dr}=-\frac{2\varepsilon}{(1+r)^2}\quad\mbox{and}\quad\xi^{\prime\prime}=\frac{d^2\xi}{dr^2}=\frac{4\varepsilon}{(1+r)^3},\]
  therefore resulting in
  \[\xi^{\prime\prime}+\frac{2}{1+r}\xi^\prime=0,\]
  on all of $\euclideanspace$ or, perhaps, we should say for all
  $r\in\{r(u)\mid u\in\euclideanspace\}=[0,\infty)$.
\end{myproof}


\bibliographystyle{plain}

\end{document}